\input amstex
\magnification=1200
\hsize=16truecm
\vsize=22truecm
\baselineskip=10pt
\centerline {\bf ON THE THETA DIVISOR OF SU(2,1).}
\hfill \par    \hfill \par
\centerline {\bf Sonia Brivio and Alessandro Verra}
\hfill \par  \hfill \par
\centerline{\bf Introduction}  \hfill \par  
Among the objects which are interestingly
related to a complex projective curve $C$,
one certainly finds the varieties
$$ SU_C(r,d).
$$  These are the moduli spaces for semistable 
vector bundles on $C$ having
rank $r$
and fixed determinant $L \in Pic^d(C)$.
In this paper we assume $C$ is smooth,
connected and of genus $g \geq 2$, then $SU_C(r,d)$ is
projective 
and its Picard group is isomorphic to $\bold Z$ (cfr. [DN]). 
Therefore we can consider the
ample generator
$H_{r,d}$ of
$Pic SU_C(r,d)$ and  the corresponding rational map
$$ h_{r,d}: SU_C(r,d) \to \bold PH^0(H_{r,d})^*.
$$ In general, as soon as $r \geq 3$, not very much is known on
the behaviour of this
map. For instance, one does'nt know under which conditions
$h_{r,d}$ is a morphism,
nor what is the degree of $h_{r,d}$ onto its image, 
(see, for a survey  on this matter, [B3]).
The situation becomes better in the case
$r = 2$. In our paper we contribute to this case
by showing the following: 
\bigskip \it { \bf THEOREM 1} $h_{2,1}$ is an embedding.
\par \rm  Of course this implies the same property for $h_{2,2d+1}$. 
The proof relies
on a method we will   explain in a moment and on a
second result we prove, which is of some independent interest:
\bigskip \it {  \bf THEOREMJ 2}JLet $F$ be a stable rank 
two vector bundle of degree $2g+1$.
If $e \in Pic^0(C)$ is general then $F(e)$ 
is globally generated and $h^0(F) = 3$.
\par \rm To put in perspective theorem 1, 
let us recall briefly some known facts
about the maps $h_{2,d}$. \par At first $h_{2,d}$
is an embedding if $d$ is even and
$C$ has general moduli.  This result is due to Laszlo, ([L1]).
If  $g$ is $\leq 3$ and $d$ is even,  $h_{2,d}$ and 
its image have been completely described by Narasimhan and
Ramanan for every curve $C$, ([NR1], [NR2]). \par If $d$ is odd and
$C$ is hyperelliptic, an explicit description of
$h_{2,d}(SU_C(2,d))$ has been given by Desale and
Ramanan, ([D-R]). From this
description it follows that $h_{2,d}$ is an 
embedding in the hyperelliptic case,
therefore the same property holds for a general $C$. 
\par Notice also that, according
to Beauville, $h_{2,d}$ has degree one unless $C$ 
is hyperelliptic of genus $\geq 3$
and
$d$ is even, (cfr. [B1] and J[B2] remarqueJ(3.14)). 
\par On the other hand, knowing
for which curves
$h_{2,d}$ fails to be an embedding is somehow a delicate question.
If $d$ is even,
a partial answer can be summarized 
in the following way (cfr. [B1], [BV]): \par
(1) Assume 
$C$ is hyperelliptic of genus $g \geq 3$,
then $h_{2,d}: SU(2,d) \to h_{2,d}(S(2,d))$
is a finite double covering. \par (2)
Assume $C$ is not hyperelliptic of genus $g
\geq 3$, then $h_{2,d}$ is injective
and its restriction to the stable locus
$SU_C^s(2,d)$ of $SU_C(2,d)$ is an embedding.
\par Since $SU_C(2,d)-SU_C^s(2,d) = Sing SU_C(2,d)$, what it is still
missing is the description of the tangent map $d(h_{r,d})$
at points of $Sing SU_C(2,0)$. \par If $d$ 
is odd the situation is less complicate and
our previous theorem 1 says that $h_{2,d}$ is always an embedding. 
Let us sketch the
method we will use for proving this theorem and the way the paper 
is organized.
\par Of course, if $d$ is odd, we can identify $SU_C(2,d)$ 
to the moduli space $X$ of
stable rank two vector bundles with determinant $\omega_C(p)$, 
where $p$ is a given
point of $C$. Then we consider in $X$ the family of open sets
$$ X_l, \quad l \in Pic^1(C),
$$ parametrizing bundles $E$ such that $F = E(l)$ 
is globally generated and has
general cohomology. One can show that for such an $F$
the natural determinant map
$$ w: \wedge^2 H^0(F) \to H^0(det F)
$$ is injective. Since $h^0(F) = 3$, the space
$$ W = Im(w)
$$ is an element of the Grassmannian $G_l$ of
$3$-dimensional subspaces of
$H^0(\omega_C(2l+p))$, ($\omega_C(2l+p) \cong det F$).
This construction defines a
rational map
$$ g_l: X \to G_l
$$ sending the moduli point of $E$ to $W$.
Since $F$ is globally generated, the
linear system $\mid W \mid$ is base-point-free. Let
$$ f_W: C \to \bold P^2
$$ be the morphism defined by $\mid W \mid$, 
it is a standard fact that
$$ f_W^*T_{\bold P^2}(-1) \cong F.
$$ Therefore we have reconstructed $F$ starting from $W$:
this implies that $g_l$
is birational and moreover that $g_l$ is biregular
at the moduli point of $E$. Assume
the two following facts are true: 
\par (a) $g_l$ is defined by a linear system which
is contained in $\mid \Cal L_X \mid$, $\Cal L_X$
being the ample generator of
$Pic X$. \par (b) $X = \cup X_l$, $l
\in Pic^1(C)$. \par 
Then we can conclude that the map associated to $\Cal L_X$ is an
embedding. \par In section 2 the map $g_l$ and its regular locus
are studied in 
detail, in particular we show (a). In section 3 we show theorem 2,
that is (b). Let
$\Theta$ be a symmetric theta-divisor in the Jacobian $J$ of $C$,
the proof of 
theorem 2 involves the pencil
$$ P_E \subset \mid 2\Theta \mid
$$ which is naturally associated to a stable vector bundle 
$E$ of determinant
$\omega_C(p)$, (see section 1 and [B2]). 
Actually theorem 2  turns out
to be equivalent to the following statement:
\par  \it for each $E$, $P_E$
has no fixed component. \rm \par \noindent
Therefore we analyze carefully the base
locus
$B_E$ of
$P_E$ and finally show that $dim B_E = g-2$.
The proof of this fact is technical.
The essential reason for it seems to be the following property:
\par \it let $Z
\subset Pic^d(C)$ be a family of subbundles of $E$
having positive degree $d$ then
$dim Z
\leq g-d-1$. \par \noindent
\rm  This property is related to Martens theorem,
(cfr. [ACGH]), and to a result
of Lange-Narasimhan saying that $E$ has
at most finitely many linear subbundles of
maximal degree, ([LN] Prop. 4.1).
\par In section 4 we study some more geometry of the base locus $B_E$,
the
main remark being that $B_E$ is reducible and 
splits in two components interchanged by
$-1$ multiplication. 
Let $W_1\subset J$ be the image of $C$ under the Abel map.
Each component of $B_E$ parametrizes
the family of translates of $W_1$
which are contained
in an element of the pencil $P_E$. In view of this,
we conclude
the paper with the following question: \par \it
describe the locus of pencils $P
\subset
\mid 2\Theta \mid$ with reducible base locus.\rm
\bigskip \noindent \it
Aknowledgments: \par \noindent
J\rm  the authors were partially supported by the
European Science Project "AGE  (Algebraic Geometry in Europe)", 
contractn. SCI-0398-C(A).
\vfill \eject

\noindent
\bf 1 Preliminary results and notations. \rm
\bigskip \noindent In this section we fix our notations.
Moreover we recollect some
results which will be frequently used.
We will always denote by
$$ X \tag 1.1
$$ the moduli space of stable rank
two vector bundles with determinant
$\omega_C(p)$, where $p$ is a given point of $C$.
$X$ is a smooth, irreducible projective variety 
of dimension $3g-3$ and $Pic X$ is isomorphic
to $\bold Z$. By definition, 
the generalized theta divisor of $X$ is the ample
generator
$$
\Cal L_X \tag 1.2
$$ of $Pic X$. The map associated to $\Cal L_X$ will be 
$$
\phi: X \to \bold P_X, \tag 1.3 
$$ where
$$
\bold P_X = \bold PH^0(\Cal L_X)^*. 
$$ By Verlinde formula $dim \bold P_X = 2^{g-1}(2^g-1)-1$. 
We reserve the notation
$$
\xi 
$$ for points of $X$. With some abuse, 
the same notation will be used for the vector
bundle corresponding to $\xi$. 
On the other hand, we will deal frequently with the
moduli space
$$ Y
$$ of semistable rank two vector bundles
of determinant $\omega_C$.  A point on $Y$,
or a bundle parametrized by this point, will be 
$$
\eta. 
$$ We recall that $Y$ is irreducible of dimension 
$3g-3$ and that $Sing Y$ is the
locus of non stable points.
$Sing Y$ is the image of the morphism $f:
Pic^{g-1}(C) \to Y$ which
sends $L \in Pic^{g-1}(C)$ to $\eta = L \oplus M$,
$M = \omega_C
\otimes L^*$. It turns out that $Sing Y$
is the Kummer variety of the Jacobian of
$C$. 
As in the case of $X$, $Pic Y$ is 
isomorphic to $\bold Z$ and its
ample generator
$$
\Cal L_Y 
$$ is called the generalized theta 
divisor of $Y$. $\Cal L_Y$ defines a map
$$
\psi: Y \to \bold P_Y,
$$ where $\bold P_Y = \bold PH^0(\Cal L_Y)^*$ 
and $dim \bold P_Y = 2^{2g}-1$. \par
Let
$$ 
\Theta \subset J = Pic^0(C) \tag 1.4
$$ 
be a symmetric theta divisor, there is a fundamental relation
between  the linear system $\mid 2\Theta \mid$
and the maps $\phi$ and $\psi$. This
relation has been explained by Beauville
in [B1] and [B2]. We need to summarize the
most important aspects of it. \par 
At first, there exists an identification
$$
\bold P_Y = \mid 2\Theta \mid \tag 1.5
$$ which relies on Wirtinger duality between
$H^0(\Cal O_J(2\Theta))$ and $H^0(\Cal
O_J(2\Theta))^*$. Let $\eta \in Y$, $\eta$ 
defines in $J$ a divisor
$$ D = \lbrace e \in J / h^0(\eta(e)) \geq 1
\rbrace, \tag 1.6
$$ with a natural structure of determinant scheme
,(cfr. [B1]). It turns out that
$$ D \in \mid 2\Theta \mid
$$ and that
$$
\psi(\eta) = D. \tag 1.8
$$ Let 
$$
\xi \in X,
$$ the next construction associates
to $\xi$ a line in $\bold P_Y$, that is a pencil
of $2\Theta$-divisors.
\bigskip \it {  \bf CONSTRUCTION 1.9} \it 
Let $\bold C_p$ be the structure sheaf of the point
$p$, $\xi$ defines the family of
exact sequences
$$
\CD 0 @>>>{ \eta_{\lambda}} @>>>
{\xi} @>{\lambda}>> {\bold C_p} @>>>{0,}\\
\endCD
$$   where $\lambda \in \bold P^1 = 
\bold PHom(\xi,\bold C_p) = \bold P{\xi^*}_p.$
\rm \par \noindent  Since $\xi$ is stable,
it follows that $\eta_{\lambda}$ is a
semistable rank two vector 
bundle of determinant
$\omega_C$. Therefore $\eta_{\lambda}$
defines the divisor
$$
\psi(\eta_{\lambda}) = D_{\lambda}
\in \mid 2\Theta \mid.
$$   We can consider the map
$$
\psi_{\xi}: \bold P^1 \to \mid 2\Theta \mid
$$ sending $\lambda$ to $D_{\lambda}$. 
One can show that $\psi_{\xi}$ is a linear
embedding, ([B2]). In particular
the line $\psi_{\xi}(\bold P^1)$ is a pencil of
$2\Theta$-divisor.
\par \rm
\bigskip \it { \bf DEFINITION 1.10}J
The previous line is the pencil
associated to $\xi$. It will be denoted by
$$ 
P_{\xi} = \lbrace D_{\lambda} /
\lambda \in \bold P^1 \rbrace.
$$
\par \rm
$P_{\xi}$ is a fundamental object for the study 
of the map $\phi$. Note that
$P_{\xi}$ is a point of the Grassmannian
$Grass (2,H^0(\Cal O_J(2\Theta))$, 
hence of its Plucker space
$$
\bold P = \bold P \wedge^2 H^0(\Cal O_J(2\Theta)). \tag 1.11
$$
\bigskip \it { \bf DEFINITION 1.12}
$ \phi_p: X \to \bold P$
is the morphism which is so defined: 
$$\phi_p(\xi) = P_{\xi}.$$
\par \rm
$\phi_p$ has been studied by Beauville in [B2]. One has
$$
\phi_p^* \Cal O_{\bold P}(1) = \Cal L_Y,
$$ therefore $\phi_p^*$ induces an
homomorphism between global sections 
$$ F_p: H^0(\Cal O_{\bold P}(1)) \to H^0(\Cal L_Y).
$$
Since the two spaces have the same dimension,
it follows  $\phi = \phi_p$ as soon as
$F_p$ is an isomorphism. It is not always true that $\phi_p$ is an
embedding nor that $F_p$ an isomorphism, (cfr. [B2]). 
A clear example is given by the
case when
$C$ is hyperelliptic and $p$ is a Weierstrass point. 
Then $\phi_p$ has degree two onto its image  
while $\phi$ is an embedding, ([B2],
3.14).
\par 
Let $\eta$ be a point of $Y$ and let 
$D = \psi(\eta)$, we want to make some remarks
on the singular locus of $D$. 
If $e \in Pic^0(C)$ we denote as
$$
b_e: C \times C \to Pic^0(C)
$$
the difference map $(x,y) \to e+x-y$.
Assume
$$
h^0(\eta(e)) \geq 2,
$$
then certainly $e \in D$.
As a corollary of Laszlo's theorem on the
singularities 
of the generalized theta divisor ([L2]),
it follows that $e \in Sing D$. 
We want to show some slight variation of this
corollary. \bigskip \noindent
{\bf PROPOSITION 1.13 }
\it Assume $b_e(C \times C)$ is not contained in $D$. 
Then:
\par \noindent
(1) $h^0(\eta(e)) = 2$
and the determinant map $v: \wedge^2 H^0(\eta(e)) \to
H^0(\omega_C(2e))$ is injective. \par \noindent
(2) Let $d = div (s)$ $=$ $\Sigma y_i$,
where $s$ is a generator of $Im (v)$. 
In $CJ\times C$
consider the divisor
$$
\Gamma = 
\Sigma C \times \lbrace y_i \rbrace,
$$ 
Then $\Gamma$ is a component of $b_e^*D$.
\par Proof. \rm (1)
Assume 
$h^0(\eta(e)) \geq 3$, then $h^0(\eta(e+x-y)) \geq 1$
and $e+x-y \in D$ for every $(x,y) \in C \times C$: 
a contradiction. Assume
$h^0(\eta(e)) = 2$ and $v$
not injective, then $\eta(e)$
contains a line bundle
$L$ with
$h^0(L) = 2$. 
Therefore $h^0(\eta(e-x)) \geq 1$
for each $x \in C$ and
hence $b_e(CJ\times C) \subset D$: 
a contradiction again. \par To show
(2) it suffices to fix a general
$x
\in C$ and compute that 
$b_e^*D$ restricted to 
$\lbrace x \rbrace \times C$
contains the restriction of
$\Gamma$. For such an $x$ 
we have $h^0(\eta(e+x)) =
2$. Otherwise it would follow
$h^0(\eta(x)) \geq 3$
generically and hence
$b_e(CJ\times C) \subset D$. 
Let $s_1,s_2$ be a basis
of $H^0(\eta(e+x))$,
it is a standard fact that 
$div (s_1 \wedge s_2)$ 
is the divisor $b_e^*D$
restricted to $\lbrace x \rbrace \times C$. 
In other words $(x-C) \cdot
D = div(s_1 \wedge s_2)$,
where $x-C$ is the image of $C$
by the Abel map
$y \to x-y$.
From the exact sequence
$$
0 \to \eta(e) \to \eta(e+x) \to \eta(e+x)_x \to 0
$$
we have $div (s_1 \wedge s_2) = d + 2x$
and the result follows. 
\bigskip
Note that $e$ is a singular point 
of the curve $b_e(\Gamma)$. Let
$$
<d> \subset \bold P^{g-1}
= \bold PH^0(\omega_C)^*
$$
be the linear span of $\phi_*d$,
where $\phi$ is the canonical map for $C$.
The projectivized tangent space
to $b_e(\Gamma)$ at $e$ contains
$<d>$, (cfr. [ACGH]). 
Since $\Cal O_C(d) = \omega_C(2e)$
it follows 
$<d> = \bold P^{g-1}$
provided $\Cal O_C(2e) \neq \Cal O_C$. 
In this case $e$ is a
singular point of every 
hypersurface $M \subset Pic^0(C)$
such that $b_e^*M$ contains
$\Gamma$. We have shown the following 
\bigskip
{\bf COROLLARY 1.14J} \it
Assume $\eta(e)$ satisfies condition (1) of 
the previous proposition. Let $M$ be a
component of
$D$ and let
$\Cal O_C(2e)
\neq
\Cal O_C$. If $b_e^*M$ contains $\Gamma$,
then $e \in Sing M$. \rm
\bigskip
We end this section
by stating some basic properties 
which are related to elementary
transformations. \par 
Let $\bold C_p$ be the structure sheaf
of a point $p \in C$, in the following 
we consider any exact sequence
$$
\CD 0 @>>> \eta_{\lambda} @>>> \xi @>{\lambda}>> 
{\bold C_p} @>>> {0,}\\
\endCD
\tag 1.15
$$ where $\eta_{\lambda}$ and $\xi$ are rank two vector
bundles. We assume that the morphism of sheaves
$\lambda$ is defined up to a non zero constant factor. 
Fixing $\xi$, the family of
all exact sequences (1.15) is parametrized by
$$
\bold P^1 = \bold PHom(\xi,\bold C_p).
$$  Let $q \subset H^0(\xi)$ be a 1-dimensional vector
space and let $q(p)$ be its
image in the fibre $\xi_p$. Assume $q(p)$ is not zero, 
then $\xi_p/q(p) \cong \bold
C_p$ and we can associate to $q$ 
the surjective morphism
$$
\pi_q = f \cdot e: \xi \to \bold C_p,
$$
where $f: \xi_p \to \xi_p/q(p)$ is the quotient map
and $e: \xi \to \xi_p$ is the
natural evaluation. $\pi_q$ is defined up 
to a non zero constant factor and its construction is linear.
That is
$$
\pi_q = \pi(q),
$$  
where
$$
\pi: \bold PH^0(\xi) \to \bold PHom(\xi,\bold C_p) \tag 1.16
$$  
is a linear map of center $\bold PH^0(\xi(-p))$. 
The Kernel of the map
$$  
h^0(\pi_q): H^0(\xi) \to H^0(\bold C_p) \tag 1.17
$$  is spanned by $q$ and $H^0(\xi(-p))$. 
\par  
The sequence (1.15) induces the
exact commutative diagram
$$
\CD 0 @>>> {H^0(\eta_{\lambda})} @>>> {H^0(\xi}) 
@>{\rho_{\lambda}}>>  {\bold C} @. \\
@. @VVV @VV{e_p}V @VVV \\ @.{(\eta_{\lambda})_p} 
@>>> {\xi_p} @>{\rho_{\lambda,p}}>>
{\bold C} @>>> 0  \\
\endCD
\tag 1.18
$$  where the vertical arrows are the evaluation maps. 
\bigskip \it { \bf PROPOSITION 1.19}\par \noindent
(1)If $\xi$ is stable, then $\eta_{\lambda}$
is semistable for each $\lambda \in \bold
P^1$.  
\par \noindent (2) If $e_p: H^0(\xi) \to \xi_p$
is surjective, then
$h^0(\eta_{\lambda})$ $=$ $h^0(\xi)-1$,
for each $\lambda \in \bold P^1$. \par
\noindent (3) Let $L \in Pic C$,
then the following conditions are equivalent:
\par
\noindent (i) $h^0(\eta_{\lambda} \otimes L) \geq 1$
for each $\lambda \in \bold P^1$,
\par
\noindent  (ii)$h^0(\xi \otimes L(-p)) = 1$
or $h^0(\xi \otimes L) \geq 2$. 
\par \rm
\it Proof \rm (1) Obvious. (2)
In the previous diagram (1.18) the
latter vertical arrow is an isomorphism.
Therefore $\rho_{\lambda}$ is surjective if
and only if $e_p$ is surjective, this implies the result.
(3) Tensor the sequence
(1.15) by $L$ and consider the corresponding diagram (1.18). 
Then
$h^0(\eta_{\lambda} \otimes L)$ $=$ $dim Ker(\rho_{\lambda}
\otimes L)$. Note that
this dimension is non zero for each $\lambda$ if and 
only if the following two conditions hold: (1) 
$h^0(\xi \otimes L)
\geq 1$, (2) either $e_pJ\otimes L$ is zero or 
$$ Im (e_p \otimes L) \cap Ker (\rho_{\lambda,p}) \neq (0) 
\quad \forall \lambda \in
\bold P^1.
$$ 
This happens if and only if $h^0(\xi \otimes L) \geq 1$
and $e_p\otimes L$ is zero
or surjective. Then the result follows.
\bigskip
\bigskip \it { \bf PROPOSITION 1.20} Assume $h^0(\xi(-y)) = 1$
for some $y \in C$. Then there
exists $\lambda \in \bold P^1 = \bold PHom(\xi,\bold C_p)$
such that the sequence
$$
\CD 0 @>>> {\eta_{\lambda}} @>>> {\xi} @>{\lambda}>>
{\bold C_p} @>>> 0 \\
\endCD
$$ is exact and $h^0(\eta_{\lambda}(-y)) \neq 0$. $\lambda$ is unique if
$h^0(\xi(-y-p)) = 0$, otherwise one has $h^0(\eta_{\lambda}(-y)) 
\neq 0$ for each $\lambda$. \par
\noindent
\par \rm
\it Proof \rm  Tensor the exact sequence (1.15)
by $\Cal O_C(-y)$ and consider
the corresponding diagram (1.18). If $e_p(-y): 
H^0(\xi(-y)) \to \xi_p(-y)$ is zero,
it follows that $h^0(\eta_{\lambda}(-y)) \neq 0$
for each $\lambda$. Moreover
$e_p(-y)$ is zero if and only if $h^0(\xi(-p-y))
\neq 0$. If $e_p(-y)$ is not zero,
then its image $q(p)$ is one dimensional. As above let 
$\pi_q: \xi(-y) \to \bold C_p$ be the surjective
morphism associated to $q(p)$.
Putting $\lambda = \pi_q(y)$ one obtains the 
unique $\lambda$ such that
$h^0(\eta_{\lambda}(-y)) \neq 0$. \bigskip 
Let
$$ 
V \subset H^0(\xi)
$$ 
be a 3-dimensional vector space and let
$$ 
w: \wedge^2 V \to H^0(det \xi) 
$$ 
be the determinant map, in the following we assume that $w$ is injective. 
\par Consider the 3-dimensional vector space
$$ W = Im(w) \subset H^0(det \xi)
$$ and the rational map
$$ f_W: C \to \bold P^2 = \bold PW^* \tag 1.21
$$ defined by $W$. Then we have
\bigskip \it { \bf PROPOSITION 1.22} Assume that
the evaluation map $V \otimes \Cal O_C \to
\xi$ is surjective, then
$$ f_W^* T_{\bold P^2}(-1) = \xi.
$$
\par \rm
\it Proof \rm  Well known. \bigskip 
\bigskip \it {\rm FURTHER CONVENTIONS AND NOTATIONS}
\rm \par \noindent
- Let $l \in Pic(C)$ and let $d \in Div(C)$. When no confusion arises  we will use
freely  $l(d)$ or $l+d$ to denote $l \otimes \Cal O_C(d)$. \par \noindent
- $C(m)$ is
the $m$-symmetric product of $C$. Let $l \in Pic^n(C)$,
then
$l-C(m)$ denotes the image of $C(m)$ in $Pic^{n-m}(C)$
by the Abel map $d \to l-d$.
\par \noindent - Let $A$ be an abelian variety and 
let $Z \subset A$. $Z_e$ is the
translate of $Z$ by $e \in A$.
$-Z = \lbrace -z, z \in Z \rbrace$. \par
\noindent -
$V^*$ is the dual of the vector bundle (vector space) $V$. 
$\bold PV$ denotes the
projectivization of $V$. $Proj(V) = \bold PV^*$. 

\par \rm
\vfill \eject
\bf 2 Very ampleness of $\Cal L_X$. \rm \bigskip 
We want to study the morphism
$$
\phi : X \to \bold P_X, \tag 2.1
$$ defined by the ample generator of $Pic X$. 
For doing this conveniently, we begin
by constructing a family of auxiliary rational maps
$$ g_l: X \to G_l, \quad l \in Pic^1(C). \tag 2.2
$$
\bigskip \it {\bf DEFINITION 2.3}JLet $(\xi,l)
\in X \times Pic^1(C)$. We say that $(\xi,l)$
satisfies condition $(+)$ if the following
three properties hold: \par \noindent
$(+1)$ $h^0(\xi(l)) = 3$, \par \noindent
$(+2)$ $\xi(l)$ is globally generated,\par \noindent
$(+3)$ the determinant map $w_l: \wedge^2 H^0(\xi(l)) 
\to H^0(det \xi(l))$ is
injective. \par \noindent By definition 
$$ X_l = \lbrace \xi \in X / (\xi,l) \quad
\text {satisfies $(+)$} \rbrace. \tag 2.4
$$
\par \rm Actually, as we are going to show,
condition $(+3)$ follows from $(+1)$
and
$(+2)$.
\bigskip \it {\bf PROPOSITION 2.5} Assume a pair $(\xi,l)$
satisfies conditions $(+1)$ and
$(+2)$, then the determinant map
$$  w_l: \wedge^2 H^0(\xi(l)) \to H^0(det \xi(l))
$$  is injective.
\par \rm
\it Proof \rm  Since $h^0(\xi(l)) = 3$, 
every vector $v \in \wedge^2 H^0(\xi(l))$ is
indecomposable. Let $v = s_1 \wedge s_2$,
$v \neq 0$. If $w_l(v) = 0$, there exists an
exact sequence
$$  0 \to A \to \xi(l) \to B \to 0 \tag 2.6
$$  such that $A$ is a line bundle and the image
of $H^0(A)$ in $H^0(\xi(l))$ contains
$s_1,s_2$. Passing to the long exact sequence,
it follows that the image of
$H^0(\xi(l))$ in $H^0(B)$ is generated by
one element $b$. Let $x$ be a
point such that $b(x) = 0$. Then the 
image of $H^0(\xi(l))$ in $\xi(l)_x$ is
contained in
$A_x$ and hence $\xi(l)$ is not globally generated. 
This contradiction shows that
$w_l$ is injective.
\bigskip 
The next proposition will be useful in various situations.
\bigskip \it {\bf PROPOSITION 2.6} For each $l \in Pic^1(C)$
there exists a projective
curve $$\Gamma_l \subset X$$ such that: \par 
(1) the set $(X-X_l) \cap \Gamma_l$ is finite, \par
(2) each point of $\Gamma_l$ satisfies conditions (+1) and (+3), \par  
(3) $h^0(\xi(l-p)) = 1$ for each $\xi \in \Gamma_l$.
\par \rm
\it Proof \rm  Observe that the set 
$$  
U_l = \lbrace M \in Pic^{g-1}(C)|
h^0(M(l-p))=0, h^0(M(l+p-x))=1, \forall x
\in C \rbrace
$$ is open and non empty for each $l \in Pic^1(C)$. Let $\eta_o = M
\oplus N$, where $M$ and $N = \omega_C \otimes M^*$ are choosen in $U$.
Then
$\eta_o(l+p)$ is globally generated, $h^0(\eta_o(l+p)) = 4$ and
$h^0(\eta_o(l-p)) = 0$. 
By semicontinuity we can replace the semistable
$\eta_o$ by a stable $\eta$
satisfying the same properties. \par
Applying construction (1.9) to $\eta(p+l)$, we obtain a family
of exact sequences
$$ 
\CD 0 @>>> {\xi_{\lambda}(l)} @>>> 
{\eta(p+l)} @>{\lambda}>> {\bold C_p} @>>>
{0,}\\
\endCD
$$ 
where $\lambda \in \bold P^1$. 
Since each $\xi_{\lambda}$ is stable, this defines
a morphism
$$
\gamma: \bold P^1 \to X
$$
such that $\gamma(\lambda) = \xi_{\lambda}$.
We denote by $\Gamma_l$ the image of
$\gamma$. \par By proposition (1.19)(2) 
we have $h^0(\xi_{\lambda}(l)) = 3$, i.e.
every point of $\Gamma_l$ satisfies property (+1). 
In principle it is possible
that $\xi_{\lambda}$
is not globally generated for every $\lambda$.
We  show that
this is not the case if $\eta$ is sufficiently general.
Let $D \in \mid 2\Theta
\mid$ be the divisor defined by $\eta$.
If $\eta$ is general $D$ is reduced and the
curve $l-C$ is transversal to $D$. Assuming this it follows that
$D \cap (l-C)$ is a set of $2g$ distinct points 
$$\lbrace l-y_1, \dots , l-y_{2g} \rbrace \quad (y_i \in C)$$
satisfying the condition $h^0(\eta(l-y_i)) = 1$. 
Moreover one has $h^0(\eta(l-x)) =
0$ for $x \neq y_1, \dots, y_{2g}$. Under these
assumptions we check the values of
$\lambda$ for which $\xi_{\lambda}(l)$ is not globally generated.
\par  If
$\xi_{\lambda}(l)$ is not globally generated at $x$, 
it follows
$h^0(\xi_{\lambda}(l-x)) = h^0(\eta(l+p-x)) = 2$.
But then $\eta(l+p-x)$ is not
globally generated at $p$ by proposition 1.19(2).
Hence $h^0(\xi_{\lambda}(l-x))=1$ and $x$
belongs to $\lbrace y_1, \dots, y_{2g} \rbrace$.
\par Since we are assuming
$h^0(\eta(l-p)) = 0$, it follows 
$h^0(\eta(l-p-x)) = 0$ for each $x \in C$. Then the evaluation  
$H^0(\eta(l-x)) \to \eta_p(l-x)$ has  one-dimensional
image $q(p)$ and   $q(p)$
defines a surjective morphism 
$\pi_q: \eta(l-x) \to \bold C_p.$, as in (1.15). 
Twisting by $\Cal O_C(p+x)$ we finally obtain
$$
\pi_q(x): \eta(l+p) \to \bold C_p.$$ 
Using the same arguments of section
1, it is easy to deduce that the
unique $\lambda$ satisfying
$h^0(\xi_{\lambda}(l-x)) = 2$ is ${\pi_q}(x)$. 
Therefore $\Gamma_l$ is a curve 
such that
$$
\Gamma_l \cap (X-X_l) = 
\lbrace \xi_{\lambda}, \lambda
= \pi_q(x), x = y_1 \dots y_{2g} \rbrace.
$$
This shows (1). 
If $\eta$ is general we can also assume that the determinant map 
$$w: \wedge^2 H^0(\eta(l+p)) \to H^0(\omega_C(2l+2p))$$
satisfies the following condition:
$Ker(w)$ does not contain indecomposable vectors.
If the genus of $C$ is $\geq 3$ 
this follows from [BV] theorem 3.3(1), where
it is shown that $Ker(w) = (0)$.
In genus two a proof can be given by a simple
dimension count, (cfr. [Br]).
The inclusion of  the subbundle $\xi_{\lambda}(l)$ in
$\eta(l+p)$ induces  the standard commutative diagram
$$
\CD
{\wedge^2 H^0(\xi_{\lambda})} @>>> {\wedge^2 H^0(\eta(l+p)} \\
@V{w_{\lambda}}VV @VwVV \\
{H^0(\omega_C(2l+p))} @>>> {H^0(\omega_C(2l+2p)),}\\
\endCD
$$
where the top arrow is injective and 
preserves indecomposable vectors. The vertical
arrows are just the determinant maps.
Since $w$ is injective, $w_{\lambda}$ must be
injective for each $\lambda$. Therefore each point of
$\Gamma_l$ satisfies property (+3) 
and the proof of (2) is completed.
To show (3)
observe that $h^0(\eta(l)) = 2$,
because we are assuming $h^0(\eta(l-p)) = 0$. Then
$\eta(l)$ is globally generated at the point $p$. From proposition (1.19)(2), it
follows $h^0(\xi_{\lambda}(l-p)) = 1$ 
for each $\lambda$, this implies (3).
\bigskip 
From (1) of the previous proposition we obtain immediately
\bigskip \it {\bf COROLLARY 2.7} For each $l \in Pic^1(C)$, $X_l$
is a non empty open set.
\par \rm
\bigskip \it {\bf REMARK 2.8} \rm For any $\xi \in X$
one can also consider the set
$$  U_{\xi} = \lbrace l \in Pic^1(C) / \text {$(\xi,l)$
satisfies $(+)$} \rbrace.
$$ It turns out that $U_{\xi}$ is not empty,
the more difficult proof of this fact
will be the main result of section 3. 
\par \rm \par Let  $l \in Pic^1(C)$, we
consider the line bundle
$$ M_l = \omega_C(2l+p) \tag 2.9
$$ and the Grassmann variety
$$ G_l = Grass(3,H^0(M_l)) \tag 2.10
$$ of $3$-dimensional subspaces. The Poincar\'e bundle
$$
\Cal E \to X_l \times C,
$$ defines on $X_l$ the sheaf
$$
\Cal H = \pi_{1*}(\Cal E \otimes \pi_2^*\Cal O_C(l)),
$$ where $\pi_i, i=1,2$ denotes
the projection onto the $i-th$ factor. Since
$h^0(\xi(l)) = 3$ for each $\xi \in X_l$, 
$\Cal H$ is a rank three vector bundle. 
The fibre of $\Cal H$ at $\xi$ is
$$
\Cal H_{\xi} = H^0(\xi(l)).
$$ A standard construction yelds a map of vector bundles
$$ w: \wedge^2 \Cal H \to H^0(M_l) \otimes \Cal O_C
$$ such that the fibrewise map
$$ w_{\xi}: \wedge^2 H^0(\xi(l)) \to H^0(M_l)
$$ is exactly the determinant map.
Note that, by our assumptions on $X_l$, $w$ is an
injective map of vector bundles. 
In particular $w$ defines a rational map
$$ g_l: X \to G_l \tag 2.11
$$ such that $$g_l(\xi) = Im(w_{\xi}),$$ 
for each $\xi \in X_l$. This is the family of
rational maps we want to consider.
\bigskip \it {\bf PROPOSITION 2.12} (1) $g_l: X \to G_l$ is birational. 
Moreover
$$ g_l/X_l: X_l \to g_l(X_l)
$$ is a biregular morphism. \par \noindent (2)
$g_l$ is defined at each point $\xi$
such that $(\xi,l)$ satisfies conditions
$(+1)$ and $(+3)$.
\par \rm
\it Proof \rm  (1) Note that $dim X = dim G_l = 3g-3$
and that both $X$ and $G_l$ are
smooth varieties. Since $X_l$ is open it suffices 
to show  that the morphism $g_l/X_l$
is injective. Then, by Zariski main theorem, 
$g_l/X_l: X_l \to g_l(X_l)$ is
biregular.  If
$\xi \in X_l$, the image of
the determinant map $w_{\xi}$ is a $3$-dimensional vector
space $W \subset H^0(M_l)$ and 
$\xi(l)$ is globally generated. The latter condition
implies that the linear system
of divisors defined by $W$ is base-point-free.
Consider the morphism
$$ f: C \to \bold P^2
$$ associated to $W$. By proposition (1.22) we have
$$ f^*T_{\bold P^2}(-1) \cong \xi.
$$ Now assume $g_l(\psi) = g_l(\xi)$ 
for some $\psi \in X_l$. Then $Im(w_{\psi}) = Im
(w_{\xi}) = W$.
Hence $\psi \cong f^*T_{\bold P^2}(-1) \cong \xi$ and $g_l/X_l$ is
injective. \par \noindent (2) 
Immediate consequence of the definitions.
\bigskip  As a next step we will
construct a commutative diagram
$$
\CD X @>{\phi_p}>> {\bold P} \\ @V{g_l}VV @VV{p_l}V  
 \\ {G_l} @>{u_l}>> {\bold P_l}  \\
\endCD
\tag 2.13
$$
such that: \par \noindent - $p_l: \bold P \to \bold P_l$
is a linear map between
projective spaces, \par
\noindent - $u_l: G_l \to \bold P_l$ is induced 
by a linear map between the Pluecker
space of $G_l$ and $\bold P_l$. \par \noindent
We recall that $\phi_p$ is the
morphism defined in 1.14 of section 1, (cfr. [B2]).
\par At first we define $p_l$: in $J$ one has the curve 
$$ C_l = \lbrace l-x, \quad x \in C \rbrace, \tag 2.14
$$ we consider the restriction map
$$
\rho_l: H^0(\Cal O_J(2\Theta)) 
\to H^0(\Cal O_{C_l}(2\Theta)) \tag 2.15
$$ and its wedge product
$$
\wedge^2(\rho_l):
\wedge^2 H^0(\Cal O_J(2\Theta)) \to \wedge^2 H^0(\Cal
O_{C_l}(2\Theta)).
$$ 
We know that the target space of $\phi_p $ is 
$$
\bold P = \bold P \wedge^2 H^0(\Cal O_J(2\Theta)),
$$ 
let
$$
\bold P_l = \bold P\wedge^2 H^0(\Cal O_{C_l}(2\Theta)), \tag 2.16
$$ 
\bigskip \it {\bf DEFINITION 2.17} $p_l: \bold P \to \bold P_l$
is the projectivization of
$\wedge^2(\rho_l)$.
\par \rm
\par  To define $u_l$ we observe that
$$
\Cal O_{C_l}(2\Theta) \cong M_l(-p), \tag 2.18
$$ leaving the proof as an exercise. 
This isomorphism yelds an obvious identity
$$
\bold P_l = \bold P\wedge^2 H^0(M_l(-p)). \tag 2.19
$$ On the other hand there exists 
a standard exact sequence
$$
\CD 0 @>>> {\wedge^3 H^0(M_l(-p))}
@>>> {\wedge^3 H^0(M_l))} @>f>> {\wedge^2
H^0(M_l(-p))} @>>> {0} \\
\endCD
$$ such that, for each indecomposable
vector $s_1\wedge s_2 \wedge s_3$, one has
$$ f(s_1\wedge s_2 \wedge s_3) = \Sigma s_i(p) 
(s_j \wedge s_k), \quad i=1,2,3.
$$
Note that $f$ is uniquely defined up to 
a non zero constant factor. Let
set
$$
\Pi_l = \bold P(\wedge^3 H^0(M_l)), \tag 2.20 
$$ $\Pi_l$ is the Plucker space of $G_l$. We 
denote the projectivization of $f$ as
$$
\overline u_l: \Pi_l \to \bold P_l. \tag 2.21
$$ Notice that
$f$ is preserving indecomposable
vectors $v = s_1\wedge s_2 \wedge s_3$. Indeed $ v = 
c(t_1 \wedge t_2 \wedge t_3)$, 
where at least two of the $t_i's$ vanish
at $p$, hence $f(v)$ is indecomposable.
Therefore the image of a point
$$ W \in G_l
$$ by $\overline u_l$ is a 
point of the Grassmannian $Grass H^0(2,M_l)$. Such a point
has an obvious geometrical interpretation, namely
$$
\overline u_l(W) = W(-p) \tag 2.22
$$ where $W(-p)$ is the subspace 
of sections of $W$ vanishing at $p$. Let
$$
\Lambda_l = \bold P\wedge^3 H^0(M_l-p)) \tag 2.23
$$ be the center of the projection $\overline u_l$, 
it is clear that
$$
\Lambda_l \cap G_l = \lbrace W \in G_l /
\text {$p$ is a base point of $\mid W \mid$
\rbrace.} \tag 2.24
$$
\bigskip \it {\bf DEFINITION 2.25} 
$u_l: G_l \to \bold P_l$ is the restriction of
$\overline u_l$ to $G_l$.
\par \rm
\par \noindent 
The next theorems 2.26 and 3.1 imply the main
theorem of this section.
\bigskip \it {\bf THEOREM 2.26} The diagram
$$
\CD X @>{\phi_p}>> {\bold P} \\ @V{g_l}VV @VV{p_l}V  
 \\ {G_l} @>{u_l}>> {\bold P_l}  \\
\endCD
$$ is commutative. 
\par \rm Let us point out
the geometric meaning of this theorem: consider the
pencil of $2\Theta$-divisors 
$P_{\xi}$ which is associated to 
$\xi \in X$. For $\xi$ general
the restriction of $P_{\xi}$ to $C_l$ is a line
$$ P  
$$ in $\mid M_l(-p) \mid$, that is a
point of the Plucker embedding of the
Grassmannian \par \noindent
$Grass(2,H^0(M_l(-p))$ in $\bold P_l$. Note that
$$ P = p_l(\phi_p (\xi)).
$$ On the other hand we consider 
the determinant map $w_l: \wedge ^2 H^0(\xi(l)) \to
H^0(M_l)$ and its image $W$. From 
$W$ we obtain the two-dimensional vector space
$W(-p) \subset H^0(M_l(-p))$. From (2.22)
we know that
$$
\bold PW(-p) = u_l(g_l (\xi)).
$$ The previous theorem says 
$$ P = \bold PW(-p). \tag 2.27
$$
\par \it Proof of theorem 2.26 \rm Let us show the 
equivalent condition (2.27). As usual
consider the family of exact sequences
$$
\CD 0 @>>> {\eta_{\lambda}(l)}
@>>> {\xi(l)} @>{\lambda}>> {\bold C_p(l)}
@>>> 0
\\
\endCD
$$ defined in (1.9). We can
assume that $\xi \in X_l$. Then, by proposition
(1.19)(2), it follows that
$h^0(\eta_{\lambda}(l)) = 2$ for each $\lambda \in
\bold P^1$. Moreover the determinant map 
$$ v_{\lambda}: \wedge^2 H^0(\eta_{\lambda}) \to H^0(\omega_C(2l))
$$ is injective because the
same holds for $w_l: \wedge^2 H^0(\xi(l)) \to
H^0(\omega_C(2l+p))$.
\par \noindent Let $\lbrace
s_{\lambda}, t_{\lambda} \rbrace$ be a basis of
$H^0(\eta_{\lambda})$ and let
$\lbrace s'_{\lambda}$, $t'_{\lambda} \rbrace $ be its
image in $H^0(\xi(l))$. Then
$$ d_{\lambda} + p = d'_{\lambda},
$$ where $d_{\lambda}$ is the divisor
of $v_{\lambda}(s_{\lambda} \wedge t_{\lambda})$
and
$d'_{\lambda}$ is the divisor
of $w_l(s'_{\lambda} \wedge t'_{\lambda})$. It
is well known that
$$ D_{\lambda}\cdot C_l = d_{\lambda},
$$ where $D_{\lambda}$ is the $2\Theta$-divisor
associated to $\eta_{\lambda}$. This
implies (2.27).
\bigskip 
\bigskip \it {\bf THEOREM 3.1} Let $\xi$ be any point of $X$. 
Then, for $l$ general in
$Pic^1(C)$, $\xi(l)$ is 
globally generated and $h^0(\xi(l)) = 3$.
\par \rm From this theorem and 
proposition (2.5) it follows that 
$$X = \cup  X_l,  \quad l \in Pic^1(C).
$$
The proof  of theorem 3.1 will be
the content of section 3.\par \noindent
Let us remind that $G_l \subset \Pi_l$
via Pluecker embedding. Therefore the
target space of $g_l$ is $\Pi_l$. We
are now in position to show that the
map 
$$ 
g_l: X \to \Pi_l
$$ 
is defined by a linear system contained 
in $\mid \Cal L_X \mid$, where $\Cal L_X$
is the generalized theta divisor.
\bigskip \it {\bf PROPOSITION 2.28} 
There exists a linear map of projective spaces
$$
\pi_l: \bold P_X \to \Pi_l
$$ such that
$$ g_l = \pi_l \cdot \phi_p.
$$
\par \rm
\it Proof \rm  The rational map 
$p_l \cdot \phi_p: X \to \bold P_l$ is defined by
$h_0, \dots , h_r$ independent sections 
of $H^0(\Cal L_X)$, ($r = dim \bold P_l$).
Let $\overline u_l: \Pi_l \to \bold P_l$ be 
the linear projection defined above in (2.21).
Since the diagram in theorem (2.26) is commutative,
we have $p_l \cdot \phi_p$ $= g_l
\cdot \overline u_l$. Therefore the map 
$g_l: X \to G_l \subset \Pi_l$ is defined by
the independent sections
$$
ch_0, \dots, ch_r, h_{r+1}, \dots h_s \in H^0(\Cal L_X^{\otimes m}),
$$
with $c \in H^0(\Cal L_X^{\otimes(m-1)})$, ($s = dim \Pi_l$). 
The linear system $\mid
H \mid$ spanned by $ch_0,\dots ,ch_r,$ $h_{r+1}, \dots ,h_s$
has no fixed components. Indeed, we
know from proposition (2.7)(2) that
$g_l$ is regular at every point $\xi$ satisfying
conditions (+1) and (+3) of definition (2.3).
Moreover, as in proposition (2.6), we
can construct a projective curve 
$\Gamma$ which is entirely contained in the set of
points satisfying (+1) and (+3).
Hence $\mid H \mid$ has no base point on
$\Gamma$. Since
$Pic X \cong \bold Z$, it follows that
$\mid H \mid$ has no fixed components. \par Now
we want to show that $c$ is constant. Let
$(x_0:\dots:x_s)$ be projective coordinates on $\Pi_l$
such that $x_0 = ch_0$,
$\dots$, $x_r = ch_r$, $x_{r+1} = h_{r+1}$, $\dots$, $x_s = h_r$
are the equations of
$g_l$. Then the center for the projection 
$\overline u_l$ is $\Lambda_l = \lbrace x_0
= \dots x_r = 0 \rbrace$. It is clear that
$$ 
g_l(div(c)) \subset \Lambda_l \cap G_l.
$$ 
On the other hand we have already remarked
in (2.24) that $W \in \Lambda_l \cap
G_l$ if and only if $p$ is a base point for
the linear system $\mid W \mid$. Therefore
every point $\xi \in div(c)$  which is not
in the base locus of $\mid H
\mid$ satisfies 
$$ 
h^0(\xi(l-p)) \geq 2.
$$ 
By semicontinuity, the latter condition holds at
each point of $div(c)$.
By proposition (2.6)(3), there exists a
projective curve $\Gamma$ 
satisfying the condition: 
$h^0(\xi(l-p)) = 1$ for each $\xi \in \Gamma$.
But then
$div(c) \cap \Gamma = \emptyset$ and $c$
is a non zero constant.
\bigskip 
Finally we can show
\bigskip \it {\bf THEOREM 2.29 } Let $C$ be any curve
of genus $g \geq 2$. Then the
generalized theta divisor of $SU(2,1)$ is very ample.
\par \rm
\it Proof \rm  We must show that $\phi$ is an embedding. 
Assume
$$
\phi(\xi) = \phi(\xi'),
$$ for two points $\xi, \xi' \in X$. By theorem
(3.1) and proposition (2.5) both $\xi$
and $\xi'$ are in the same open set $X_l$,
provided $l$ is sufficiently general.
Since, by proposition (2.12), $g_l/X_l$ 
is injective it follows $\xi \cong \xi'$.
Hence
$\phi$ is injective. \par
\noindent Assume
$$ (d\phi)_{\xi}(v) = 0
$$ for a point $\xi \in X$ and a tangent
vector $v \in T_{X,\xi}$. Note that the
linear projection $\pi_l$ is defined at 
$\phi(\xi)$, provided $l$ is general. Indeed
$\pi_l$ is defined at $\phi(\xi)$ if no 
element of the pencil $P_{\xi}$ contains
$C_l$. The latter condition is true for
a general $l$. This is an immediate consequence of 
the following property: 
the family of curves $C-l$ which 
are contained in an element $D \in P_{\xi}$ has
dimension $\leq g-2$, (cfr. [BV] 5.10).
Therefore we have
$$ (d\pi_l)_{\phi(\xi)} \cdot
(d\phi)_{\xi} = (dg_l)_{\xi}.
$$ By theorem (3.1) $\xi \in
X_l$ for $l$ general. Then, by proposition (2.12)(1),
$(dg_l)_{\xi}$ is injective.
Hence $v=0$ and $(d\phi)_{\xi}$ is injective. 
\bigskip 

\vfill \eject
\bf 3 The main technical result. \rm \bigskip
The purpose of this section is to show
the following theorem:
\bigskip \it {\bf THEOREM 3.1}JLet $\xi$ be a stable rank
two vector bundle on C and let
$det \xi \in Pic^{2g-1}(C)$. Then, for a general $l \in Pic^1(C)$,
the following
conditions are satisfied: \par \noindent - $h^0(\xi(l)) = 3$,
\par \noindent -
$\xi(l)$ is globally generated. 
\par \rm It suffices to show the theorem 
when $det \xi = \omega_C(p)$, ($p \in
C$), therefore we will assume as usual 
$$
\xi \in X.
$$ The proof of the theorem is given at
the end of the section. Preliminarily, we
introduce some definitions and lemmas.
\bigskip \it {\bf DEFINITION 3.2}JLet $\xi \in X$, we
fix the following notations: \bigskip
\noindent  - $E_{\xi} = \lbrace e \in Pic^0(C) / h^0(\xi(e)) \geq 2 \rbrace. $
\bigskip
\noindent  - $ V_{\xi} =
\lbrace e \in Pic^0(C) / h^0(\xi(e-y)) \geq 1,
\forall y \in C \rbrace.$J\bigskip \noindent
- $H_{\xi} = \lbrace l \in Pic^1(C) /
l-C \subset E_{\xi} \rbrace.$
\bigskip
\noindent By definition $E_{\xi}$ is the exceptional locus of $\xi$.
\par \rm Let $e \in E_{\xi}$.
Tensoring by $\Cal O_C(e)$ the usual sequence 
$$
\CD 0 @>>> {\eta_{\lambda}} @>>> {\xi} @>{\lambda}>> {\bold C_p} @>>> 0 @. {\quad
\quad (\lambda \in \bold P^1)} \\ 
\endCD
\tag 3.3
$$ and passing to the long exact sequence, we obtain
$$ 0 \to H^0(\eta_{\lambda}(e)) \to H^0(\xi(e)) \to \bold C \to ... .
$$ Then it follows $h^0(\eta_{\lambda}(e)) \geq 1$ 
and hence $e \in D_{\lambda}$,
where $D_{\lambda}$ is the divisor of $\eta_{\lambda}$.
We have shown the following
\bigskip \it {\bf LEMMA 3.4}\par \noindent (i) 
$E_{\xi}$ is in the base locus of the pencil
$P_{\xi}$.
\par \noindent (ii) $dim E_{\xi} \leq g-1$.
\par \rm 
Our main problem will be to exclude $dim E_{\xi} = g-1$.
\bigskip \it 
{\bf REMARK}J\rm  Let $\Cal P \to Pic^0(C) \times C$
be a Poincar\'e bundle and
let $\Cal E =
\Cal P \otimes \pi_2^*\xi$, where $\pi_i$ 
denotes the projection onto the
$i$-th factor. Of course we can consider 
the standard exact sequence
$$ O \to \Cal E(-D) \to \Cal E 
\to \Cal E \otimes \Cal O_D \to 0,
$$ where $D = Pic^0(C) \times \lbrace p \rbrace$. 
Applying the functor $\pi_{1*}$ to
this sequence, it follows that $E_{\xi}$ is 
the Support of $R^1\pi_{1*}\Cal E$, in
particular $E_{\xi}$ has dimension $\geq g-2$.
We will study $E_{\xi}$ in the
next section. 
\par \rm
\bigskip \it 
{\bf LEMMA 3.5} $dim H_{\xi} \leq g-2$. 
\par \rm
\it
Proof \rm  Let $l \in H_{\xi}$ and let $y \in C$. 
Tensoring  (3.3) by $\Cal
O_C(l-y)$ and passing to the long exact sequence, we obtain
$$  0 \to H^0(\eta_{\lambda}(l-y))
\to H^0(\xi(l-y)) \to \bold C \to ... .
$$ Then $h^0(\eta_{\lambda}(l-y)) \geq 1$
and the curve $l-C$ is contained in the
divisor $D_{\lambda}$ of $\eta_{\lambda}$. 
This shows that $H_{\xi} \subset
H_{\eta_{\lambda}} = \lbrace l \in Pic^1(C)
/ l-C \subset D_{\lambda}
\rbrace$. On the other hand it is known
that $dim H_{\eta_{\lambda}} \leq g-2$, ([BV]
5.10). This implies the result. 
\bigskip 
\bigskip \it 
{\bf COROLLARY 3.6} $h^0(\xi(l))$ $=$ $3$
for a general $l \in Pic^1(C)$.
\par \rm
\it Proof \rm  Let $H^1_{\xi} = \lbrace l \in Pic^1(C) / h^0(\xi(l))
\geq 4 \rbrace$.
For each $l \in H^1_{\xi}$ the curve $l-C$
is contained in $E_{\xi}$. Hence
$H^1_{\xi}$ is contained in $H_{\xi}$ and,
by the previous lemma, $dim H^1_{\xi} \leq
g-2$.
\bigskip 
\bigskip \it {\bf LEMMA 3.7}
\par \noindent (i)
$V_{\xi} \subset E_{\xi}$, \par
\noindent  (ii) $dim V_{\xi} \leq g-2$. \par \rm 
\it Proof \rm  (i) By definition a point
$e$ in $V_{\xi}$ satisfies the condition
$h^0(\xi(e-y)) \geq 1$, $\forall y \in C$. 
It si obvious that this implies
$h^0(\xi(e)) \geq 2$.
\par \noindent (ii) Let $V \subset V_{\xi}$
be an irreducible component of dimension
$\geq g-1$, consider the difference map 
$$
\delta : V \times C \to Pic^{-1}(C), \quad (\delta(e,y) = e-y).
$$ Since $dim V \geq g-1$ it
follows that $\delta$ is surjective. The surjectivity of
$\delta$ implies
$$ h^0(\xi(-l)) \geq 1 \quad \forall l \in Pic^1(C). 
$$ From Riemann-Roch,  Serre 
duality and the isomorphism $\omega_C \otimes \xi^*(l)
\cong \xi(l-p)$ it follows that
the previous condition is equivalent to
$$ h^0(\xi(l-p)) \geq 2 \quad \forall l \in Pic^1(C).
$$ But then $E_{\xi} = Pic^1(C)$, against lemma (3.4). 
\bigskip  In the next proposition we 
bound the dimension of a family of positive linear
subbundles of $\xi$.
\bigskip \it {\bf PROPOSITION 3.10}  Let
$Z \subset Pic^d(C)$ be a closed subset and
let $d
\geq 1$. Assume 
$$   h^0(\xi \otimes L^*) \geq 1,
\forall L \in Z. \tag 3.11 $$ 
Then $d \leq g-1$ and moreover $dim Z \leq g-1-d$.
\par \rm
\it Proof \rm  We will say that any 
closed set satisfying (3.11) is a family of
subbundles of $\xi$. By stability $d \leq g-1$,
let us show that $dim Z \leq g-d-1$.
\par
\noindent   There is no restriction to assume $Z$ irreducible. 
To obtain a
contradiction we assume  
$$ dim Z \geq g-d,
$$  then we consider the closed set
$$
\tilde Z = \lbrace (L,x) \in Z \times
C / h^0(\xi\otimes L^*(-x)) \geq 1 \rbrace
$$  and its first projection
$$  Z^n \subset Z.
$$ 
\bigskip \it {\bf CLAIM} There exists a family
$T \subset Pic^k(C)$ of subbundles of $\xi$
such that 
$$ dim T \geq g-k \quad \text {and} \quad T-T_n \neq \emptyset. $$ 
Therefore, up to
replacing $Z$ by $T$, we can assume 
$$  Z-Z^n \neq \emptyset.
$$
\par \rm \noindent To show the claim assume $Z = Z^n$
and consider the sum map
$$
\phi: \tilde Z \to Pic^{d+1}(C), \quad (\phi(L,x) = L(x)).
$$
$\phi$ defines a family
$$  Z_1 = \phi(\tilde Z)
$$  of subbundles of $\xi$ of degree $d+1$.
Since $Z = Z^n$, it follows
$dim \tilde Z$ $\geq dim Z$.
On the other hand each fibre of $\phi$ is at most one
dimensional, hence
$$  dim Z_1 \geq dim \tilde Z -1 \geq dim Z -1 \geq g-(d+1).
$$  If $Z_1-Z_1^n$ is non empty we replace $Z$
by the family $T = Z_1$. If $Z_1-Z_1^n$
is empty we go on with the same construction.
After a finite number $s$ of steps, we
obtain the required family $T$
in $Pic^{d+s}(C)$. This shows the claim.
\par \noindent 
Let $e\in Pic^0(C)$ and let $Z_e$  be the translate of $Z$ by $e$.
Then $Z_e$ is a
family of subbundles of $\xi(e)$. By lemma 3.4 $h^0(\xi(e)) = 1$ if
$e$ is general. Therefore, up to a translation, we can assume 
$$   h^0(\xi) = 1.
$$   Let
$$   W_d
$$   be the image of the Abel map $a: C(d) \to Pic^d(C)$.
Up to translating $Z$ by $e$,
we can also assume that $W_d$ is transversal
to $Z$ and $Z_n$. As a consequence of the
previous assumptions it follows that
$$  W_d \cap (Z-Z_n)
$$  is non empty. We want to remark 
that $W_d \cap (Z-Z_n)$ contains two distinct
points. This is obvious if $dim Z > g-d$. 
Assume $dim Z = g-d$. Then, by
transversality, we have $W_d \cap Z^n = \emptyset$. 
Moreover the cardinality of $W_d
\cap Z$ equals the intersection index $(W_d,Z)$. 
It is well known that the latter
cannot be one in the Jacobian of $C$.
\par \noindent   Now we can complete the proof: let
$A_1,A_2$ be distinct points of
$W_d \cap (Z-Z^n)$, then $A_i = \Cal O_C(a_i)$,
where $a_i$ is effective. We have
$h^0(\xi(-a_i)) = h^0(\xi) = 1$, $i=1,2$. 
Hence a non zero section $s \in H^0(\xi)$
vanishes on $Supp a_1 \cup Supp a_2$. 
But then $h^0(\xi(-a_1-x)) = 1$ for some $x \in
Supp a_2$ and $A_1 \in Z^n$: a contradiction.
\bigskip 
\bigskip \it {\bf REMARK}\rm  
 The proposition implies 
the following property: Let $E$ be a stable rank
two vector bundle
and let $s(E) = deg E - 2d$, where $d$ is
the maximal degree for a linear line subbundle of $E$. 
If $s(E) = 1$ the family $Z$ of all subbundles
of degree $d$ is finite. Indeed, tensoring $E$
and its maximal subbundles by a suitable $L \in
Pic(C)$,  we can assume $det E = \omega_C(p)$
and $d = g-1$. The property was already proved by
Lange-Narasimhan in [LN] 4.2, where they also compute 
the number of maximal subbundles.
\par \rm Now we fix an irreducible component
$$ B
$$ of the exceptional locus $E_{\xi}$
and we assume that
$$ dim B = g-1. \tag 3.12
$$
\bigskip \it {\bf LEMMA 3.13} If $B$
exists a general point $e \in B$ satisfies the
following conditions (*) and (**).
\par (*) $h^0(\xi(e)) = 2$ and the determinant map 
$$  v: \wedge^2 H^0(\xi(e)) \to H^0(det \xi(e))
$$ is injective, \par (**) let $d_e \in \mid det \xi(e) \mid$
be the
divisor defined by $Im(v)$. 
For each $y \in Supp d_e$ the curve
$$ C_y = \lbrace e+x-y / x \in C \rbrace
$$ is not contained in $B$. 
\par \rm
\it Proof \rm  Let $e \in Pic^0(C)$,
it is easy to check that $e \in V_{\xi}$ if and
only if one of the following
conditions is satisfied: \par \noindent  (1)
$h^0(\xi(e)) \geq 3$,
\par \noindent  (2) $h^0(\xi(e)) = 2$
and the determinant map $v:J\wedge^2 H^0(\xi(e))
\to H^0(det \xi(e))$ is zero. \par \noindent
Assume $e$ is a general point of $B$.
Then $e$ is not in $V_{\xi}$ because $dim V_{\xi} \leq g-2$. 
Hence $h^0(\xi(e)) = 2$
and the map
$$ v: \wedge^2 H^0(\xi(e)) \to H^0(det \xi(e))
$$ is injective. This shows that $e$ satisfes (*).  \par 
To show that $e$
satisfies (**), consider the two natural projections
$$
\CD C @<{\alpha}<< {C \times C}J@>{\beta}>> C \\
\endCD
\tag 3.14
$$ and the exact sequence
$$ 0 \to \alpha^*\xi(e-y) \to
\alpha^*\xi(e-y) \otimes \Cal O_{C \times C}(\Delta)
\to \alpha^*\xi(e-y) \otimes \Cal O_{\Delta}(\Delta) \to 0,
\tag 3.15
$$ where $\Delta \subset C \times C$ 
is the diagonal and $y \in Supp d_e$. \par
Applying the functor $\beta_{*}$
to the previous sequence, one obtains:
$$ 0 \to H^0(\xi(e-y))\otimes \Cal O_C \to R^0 
\to \omega_C^*\otimes \xi(e-y) \to
H^1(\xi(e-y)) \otimes \Cal O_C
\to R^1 \to 0,
\tag 3.16
$$ where
$$ R^i =  R^i\beta_*(\alpha^*\xi(e-y) 
\otimes \Cal O_{C \times C}(\Delta)).
$$ Since $y \in Supp d_e$, we have $h^0(\xi(e-y)) \geq 1$
and hence
$h^1(\xi(e-y)) \geq 2$. Moreover, for a general 
$x \in C$, it holds
$$ R^0_x = H^0\xi(e+x-y)) \quad , 
\quad R^1_x = H^1(\xi(e+x-y)).
$$ To obtain a contradiction let's assume 
$$ C_y \subset B
$$ for a general point $e \in B$.
Then the rank of $R^0$ is $2$ and  the rank of
$R^1$ is $1$. This implies that, 
in the previous sequence 3.16, the image of  the
coboundary map 
$$ c: \omega_C^*\otimes \xi(e-y) 
\to H^1(\xi(e-y)) \otimes \Cal O_C
$$ has rank (at most) one.
On the other hand, by Serre 
duality, the dual map 
$$ c^*: H^0(\xi(e+y-p)) \otimes \Cal O_C \to \xi(e+y-p).
$$ is just the evaluation map. Since the sheaf $Im(c) $
has rank one $c^*$ is not
generically surjective. Therefore the 
image of $c^*$ is a line bundle
$$ N \subset \xi(e+y-p) \tag 3.17
$$ such that $h^0(N) \geq 2$. In particular it follows that
$$ e+y-p \in V_{\xi}.
$$ We want to show that this yelds a contradiction.
Consider the closed sets
$$
\tilde B = \lbrace (e,y) \in B \times C / y \in Supp d_e \rbrace
$$ and 
$$ T = \lbrace e+y-p \in Pic^0(C) / (e,y) \in \tilde B \rbrace .
$$ From the previous remarks it is clear that $T \subset V_{\xi}$.
On the other hand 
$$ T =\phi(\tilde B),
$$ where $\phi: \tilde B \to Pic^0(C)$ is the morphism
sending $(e,y)$ to $e+y-p$. The
fibre of $\phi$ at $f = \phi(e,y)$ is 
$\phi^{-1}(f) = \lbrace (e+y'-y,y'), y'J\in C
\rbrace$. 
Then, since $dim \tilde B \geq g-1$, it follows 
$$ dim T = g-2.
$$ Let us see that this is impossible. Consider the closed set
$$ W = \lbrace (N,f) / f \in T, NJ\subset \xi(f),
N \in W^1_d \rbrace, 
$$ where $W^1_d \subset Pic^d(C)$, $2 \leq d \leq g-1$,
is the Brill-Noether locus.
From the previous remarks 
it follows that the projection of $W$ onto $T$
is surjective. Hence
$$ dim W \geq g-2.
$$ The difference map
$$
\delta: W \to Pic^d(C)
$$ yelds a family 
$$ Z = \delta(W)
$$  of linear subbundles of $\xi$ of positive degree.
For the fibre of $\delta$ at a
point $L = \delta(N,f)$ we have
$$
\delta^{-1}(L) \subseteq \lbrace (N',N-N'+f) /
N' \in W^1_d, \rbrace \cong W^1_d.
$$ 
By Martens theorem $dim W^1_d \leq d-2$, 
($dim W^1_d \leq d-3$ if $C$ is not
hyperelliptic). This implies $dim Z \geq g - d$,
against proposition 3.10. 
\bigskip 
\bigskip \it {\bf LEMMA 3.18} If $B$ exists
a general point
$e \in B$ satisfies  the
following condition (***): \par \noindent
(1) There exist two points $y_1,y_2 \in C$
such that
$$
h^0(\xi(e-y_1)) = h^0(\xi(e-y_2)) = 1 \quad , 
\quad h^0(\xi(e-y_1-y_2)) = 0.
$$
(2) If $h^0(\xi(e-p)) = 0$, there exist two distinct points 
$\lambda_1, \lambda_2 \in \bold
PHom(\xi,\bold C_p)$ such that 
$$ h^0(\eta_{\lambda_i}(e-y_i)) = 1 \quad i=1,2
$$
where $\eta_{\lambda_i}(e)$ is defined by $\lambda_i$ via the 
exact sequence
$$
\CD 0 @>>> {\eta_{\lambda_i}(e)} @>>> 
{\xi(e)} @>{\lambda_i}>>{\bold C_p(e)} @>>>
{0.}\\
\endCD
$$ 
\par \noindent
(3) If $h^0(\xi(e-p)) \geq 1$, there exists $\lambda \in 
\bold PHom(\xi,\bold C_p)$ such that $h^0(\eta_{\lambda}(e)) = 2$,
where $\eta_{\lambda}(e)$ is defined by $\lambda$ via the  exact sequence
$$
\CD 0 @>>> {\eta_{\lambda}(e)} @>>>  {\xi(e)} @>{\lambda}>>{\bold C_p(e)} @>>>
{0.}\\
\endCD
$$
In particular the determinant map $v_{\lambda}: \wedge^2 H^0(\eta_{\lambda}(e)) \to
H^0(\omega_C(2e))$ is injective
\it Proof \rm  By lemma 3.13 conditions (*) and (**)
are satisfied on a dense open set
$$ 
U \subset B.
$$ 
For $e \in U$  the determinant 
map $v_e: \wedge^2 H^0(\xi(e))
\to h^0(det \xi(e))$ is injective and $h^0(\xi(e))$ $= 2$.
Therefore $e$ defines
the effective divisor
$$ 
d_e = div(s_e) \in \mid \omega_C(2e+p) \mid, 
$$
where $s_e$ is a generator of $Im(v_e)$. \par
{\bf CLAIM 3.19} \it Assume $e$ is general in $U$. Then
$Supp (d_e)$ contains at least $g-1$ distinct points
$y_1, \dots y_{g-1}$, these points are different from $p$. \rm
\par
\it Proof of the claim. \rm
We consider the morphism $b: U \to C(2g-1)$
sending $e$ to $d_e$ and the Abel map  
$a: C(2g-1) \to Pic^{2g-1}(C)$. Note that  
$$ 
a \cdot b: U \to Pic^{2g-1}(C)
$$ 
is simply the restriction to $U$ of the morphism
$$ 
h: Pic^0(C) \to Pic^{2g-1}(C),
$$ 
where $h(e) = det\xi(e) = \omega_C(2e+p)$. $h$ is finite
of degree $2^{2g}$,  therefore
$$ 
dim  b(U) = dim U = g-1.
$$ 
Let 
$$ d_e = k_1y_1 + \dots + k_ry_r, 
$$ 
where $y_1, \dots y_r$ are distinct points
and $k_i \geq 1$. Then, for $e$ general in $U$, 
$r \geq g-1$:
otherwise we would have $dim$ $b(U) < g-1$. For the same
reason $r \geq g$ if $Supp d_e$ constantly contains
a point $p$. \par
Let $y_1 \in Supp d_e$, if $h^0(\xi(e-y_1))
\geq 2$ the curve 
$$
C_{y} = \lbrace e+x-y_1, x \in C \rbrace
$$
is contained in $E_{\xi}$. Since $e$ is general
we can assume that $B$ is the only
component of $E_{\xi}$ through $e$. Hence
$C_{y} \subset B$, which is impossible because
$e$ satisfies condition (**) of lemma (3.13). 
This contradiction implies that $h^0(\xi(e-y_1)) = 1$,
($\forall y_1 \in Supp d_e$).
Now assume $d_e =
k_1y_1 + \dots + k_ry_r$ with $r > g-1$,
then fix a non zero section
$s_1$ vanishing at $y_1$. $s_1$ defines a
saturated subbundle $\Cal O_C(d_1)$ of $\xi(e)$, where $d_1$ 
is an effective divisor
such that $y_1 \leq d_1 \leq d_e$. 
By stability $deg d_1 \leq g-1$. Hence
the set $Z = Supp d_e - Supp d_1$ is not empty and 
we can  choose a non zero section $s_2$ vanishing at
a point $y_2 \in Z$. By construction
$\lbrace s_1,s_2 \rbrace$
is a basis of $H^0(\xi(e))$, moreover $s_1(y_2)$,
$s_2(y_1)$ are not zero. Therefore
$h^0(\xi(e-y_1-y_2)) = 0$. 
\par Finally we assume $r = g-1$.
Then, for general $e$, $Supp d_e$ consists of
$g-1$ distinct  points. 
Fixing $s_1$ and $d_1$ as above, we have that either $d_1 =
x_1 +
\dots x_{g-1}$ or $Supp d_e - Supp d_1 \neq \emptyset$.
In the latter case we repeat the previous argument.
In the former case we can construct as in
section 1 an exact sequence
$$ 
0 \to \eta(e-d_1) \to \xi(e-d_1) \to \bold C_p(e-d_1) \to 0
$$ 
such that $h^0(\eta(e-d_1)) \geq 1$.
Then $\eta$ is semistable
not stable and its $2\Theta$-divisor is
$D = \Theta_a + \Theta_{-a}$,
for some $a \in Pic^0(C)$. Since $D \in P_{\xi}$ and no
component of $D$ can move in a pencil,
it follows that the base locus of $P_{\xi}$ has
codimension two. Hence $B$ cannot exist and this case is
impossible. This shows part (1) of condition (***). 
To show part (2) we take the two previous sections 
$s_1,s_2 \in H^0(\xi(e))$. Since $s_i$ vanishes at $y_i$,
there exists an exact sequence
$$
\CD 0 @>>> {\eta_{\lambda_i}(e)} @>>> 
{\xi(e)} @>{\lambda_i}>>{\bold C_p(e)} @>>>
{0,}\\
\endCD
$$ 
such that $h^0(\eta_{\lambda_i}(e-y_i)) \geq 1$, ($i=1,2)$.
This sequence is constructed from the evaluation of $s_i$
at $p$ as in section 1, (1.20). Since $h^0(\xi(e-p)) = 0$,
we have $\lambda_1 \neq \lambda_2$ and 
$h^0(\eta_{\lambda_i}(e-y_i)) = 1$.
This follows from (1.16) and proposition 1.19(2). \par
To show part (3) observe that the evaluation map
$e_p: H^0(\xi(e)) \to \xi(e)_p$ is one dimensional 
because $h^0(\xi(e-p)) = 1$. Then the image of $e_p$
defines the required $\lambda$ as in (1.18). In particular, it is obvious
that the determinant $v_{\lambda}$ is injective, because the same holds
for $v: \wedge^2 H^0(\xi(e)) \to H^0(\omega_C(2e+p))$.
\bigskip \it {\bf THEOREM 3.20} 
$dim E_{\xi} \leq g-2$ for each $\xi \in X$. \par \rm
\it Proof \rm  
Let $B \subset E_{\xi}$ be an irreducible component of dimension $g-1$
and let
$$ 
P_{\xi} = \lbrace D_{\lambda}, \quad \lambda \in \bold P^1 \rbrace
$$ 
be the pencil associated to $\xi$. 
Since $E_{\xi}$ is in the base locus of
$P_{\xi}$, we have 
$$  
D_{\lambda} = M_{\lambda} +
F \quad \forall \lambda \in \bold P^1,
$$  
where $\lbrace M_{\lambda}, \lambda \in \bold P^1 \rbrace$ 
is a pencil with no
fixed components and $F$ is an effective 
divisor containing $B$.
\par Let $e$ be a general point of $B$, 
we can assume that $e$ satisfies the
conditions (*), (**), (***) 
which have been defined in lemmas 3.13 and 3.18. Moreover
we can also assume that $e$
satisfies the following conditions:
\par \noindent  
- $e$ is not in the base locus of the pencil $P_M = \lbrace M_{\lambda},
\lambda \in \bold P^1J\rbrace$, \par \noindent
- $e$ is smooth for the unique element of $P_M$ passing
through $e$, \par \noindent
- $B$ is the unique irreducible component of $F$
passing through $e$.
\par
Assume $h^0(\xi(e-p)) = 0$. Since (***)(2) holds, there are
two distinct points 
$\lambda_1, \lambda_2 \in \bold PHom(\xi,\bold C_p)$ such that 
$$ h^0(\eta_{\lambda_i}(e-y_i)) = 1 \quad i=1,2
$$ where $\eta_{\lambda_i}(e)$ is defined by $\lambda_i$ 
via the  exact sequence
$$
\CD 0 @>>> {\eta_{\lambda_i}(e)} @>>>
{\xi(e)} @>{\lambda_i}>>{\bold C_p(e)} @>>>
{0.}\\
\endCD
$$
Now observe that the condition $ h^0(\eta_{\lambda_i}(e-y_i)) = 1$ implies 

$$
C_{y_i} = e+C-y_i \subset D_{\lambda_i}, \quad (i=1,2).
$$ 
Note that the curve $C_{y_i}$ is not contained in
$F$. Indeed $B$ does not contain $C_{y_i}$
because $e$ satisfies (**). On the other
hand $e$ is a point of $C_{y_i}$
and every irreducible component
of $F$ which is different from $B$
does not contain $e$. Therefore $C_{y_i}$ is not
contained in $F$ and hence 
$$
C_{y_i} \subset M_{\lambda_i}, \quad (i=1,2).
$$ 
This implies that $e$ is in the base
locus of the pencil $\lbrace M_{\lambda}, \lambda
\in \bold P^1 \rbrace$: a contradiction. \par
Finally, assume $h^0(\xi(e-p)) \geq 1$, then we know from 
condition (***)(3) that there exists an exact sequence
$$
\CD 0 @>>> {\eta_{\lambda}(e)} @>>>  
{\xi(e)} @>{\lambda}>>{\bold C_p(e)} @>>> {0.}\\
\endCD
$$ 
such that $h^0(\eta_{\lambda}(e)) = 2$
and the determinant map $v_{\lambda}$
is injective. Obviously
we can choose $\Cal O_C(2e) \neq \Cal O_C$,
so that $\eta_{\lambda}(e)$ satisfies
all the assumptions of corollary
1.14. 
Let $M = M_{\lambda}$ be the corresponding
element of the pencil
$P_M$. 
With the same notations of Proposition (1.13),
we consider 
the curve
$$ b_e^*(M+F) \subset C \times C.$$
Recall that $b_e^*(M+F)$ contains
the divisor  $\Gamma = \Sigma C \times
\lbrace y_i \rbrace$; where, in
our case, $\Sigma y_i = d_e-p$. 
Since $F$ does not contain
$C_{y_i}$, $b_e^*F$ does not contain
$C \times \lbrace y_i \rbrace$. Hence 
$\Gamma$ is contained in
$b_e^*M$. Then, by the corollary (1.14),
$e$ is a singular point of $M$: a
contradiction.

\bigskip  Finally we can give a
\bigskip \it {\bf PROOF OF THEOREM 3.1}J\par \rm
\it Proof \rm  
Let $\sigma: Pic^0(C) \times C \to Pic^1(C)$ be the sum map. 
It is clear
that 
$$
\sigma(E_{\xi} \times C) = Y, 
$$ where
$$ Y = \lbrace l \in Pic^1(C)/ h^0(\xi(l-x)) 
\geq 2 \text {, for some $x \in C$}
\rbrace.
$$ On the other hand $l \in Y$ 
if and only if $\xi(l)$ is not globally generated or
$h^0(\xi(l)) \geq 4$. By the previous 
theorem 3.20 we have $dim Y \leq dim(E_{\xi}
\times C) \leq g-1$. 
Hence the complement of $Y$ is not empty
and theorem 3.1 follows.
\bigskip 
\vfill \eject
\bf 4 Pencils of $2\Theta$-divisors. \rm 
\bigskip \noindent In this section we
describe some geometry of the base locus
$$ B_{\xi} \tag 4.1
$$ of a pencil $P_{\xi}$. 
In particular we will see that $B_{\xi}$ is reducible if
$\xi$ is general,  
a component of it being the exceptional locus 
$$E_{\xi} = \lbrace e \in Pic^0(C) / h^0(\xi(e)) \geq 2
\rbrace. \tag 4.2 $$
$B_{\xi}$ is considered with its natural structure of scheme,  
it is useful to recall
that
$B_{\xi} = -B_{\xi}$.
\bigskip \it {\bf LEMMA 4.3} For each $\xi$
one has 
$$ 
B_{\xi} = E_{\xi} \cup -E_{\xi}.
$$
In particular $E_{\xi}$ is a component of 
$B_{\xi}$, moreover
$E_{\xi} \neq -E_{\xi}$ if $\xi$
is general.
\par \rm
\it Proof \rm It is clear that $E_{\xi} \cup -E_{\xi}
\subset B_{\xi}$. Let $e \in B_{\xi}$
be in the complement of $E_{\xi}$, 
then $h^0(\xi(e)) = 1$. Since $e$ is 
in the base locus of $P_{\xi}$,
we have $h^0(\eta_{\lambda}(e)) = 1$ for each $\lambda
\in \bold P^1$. Hence,
from proposition (1.19), it follows $h^0(\xi(e-p)) = 1$
and finally $h^1(\xi(e-p)) = 2$.
But then $e \in -E_{\xi}$ because 
$h^0(\xi(-e))$ $=$ $h^1(\xi(e-p)) = 2$. 
Hence $B_{\xi} = E_{\xi} \cup -E_{\xi}$.
\par \noindent 
To complete the proof, let us produce one $\xi$ such that 
$E_{\xi} \neq -E_{\xi}$. 
At first there is no
problem to construct a stable, globally generated rank 
two vector bundle $\eta(l)$
having very ample determinant $\omega_C(2l) \in Pic^{2g}(C)$ such 
that $h^0(\eta(l))= 3$. 
For this just
choose $\eta(l) = f_W^*T_{\bold P^2}(-1)$, where 
$f_W: C \to \bold P^2$
is the
morphism associated to a general $3$-dimensional subspace 
$W \subset
H^0(\omega_C(2l))$. We can also assume that $f_W(C)$ 
has no cusp so
that $h^0(\eta(l-2p)) = 0$. Let $e = l-p$,  
by (1.19) $\eta(l)$ induces a family of exact sequences
$$
\CD 0 @>>> {\xi_{\lambda}(e)} @>>>
{\eta(l)} @>{\lambda}>> {\bold C_p} @>>> {0,}\\
\endCD
$$
where $\xi_{\lambda}$ is stable and 
$h^0(\xi_{\lambda}(e)) = 2,$  $\forall
\lambda \in \bold P^1$. Passing to the corresponding
long exact sequences we have 
$h^1(\xi_{\lambda}(e)) = 1$
iff $h^0(\lambda): H^0(\eta(l)) \to \bold C_p$
is  surjective. On the other hand we have 
$h^0(\xi_{\lambda}(-e)) = h^1(\xi_{\lambda}(e))$.
Note that $h^0(\eta(e)) = 1$
and that $h^0(\eta(e-p)) = 0$. Therefore,
applying proposition (1.20), there
exists a unique $\lambda_o$
such that $h^0(\lambda_o)$ is not surjective.
Let $\lambda \neq \lambda_o$, then
$h^0(\xi_{\lambda}(-e)) = 1$ and hence  
 $E_{\xi_{\lambda}} \neq -E_{\xi_{\lambda}}$.
\bigskip 
\bigskip \it {\bf REMARK}J\rm It is possible that 
$E_{\xi} = -E_{\xi}$ in some special situations.
For instance if
$C$ is hyperelliptic and
$\xi \cong i^*\xi$, ($i$ hyperelliptic involution),
it follows $-E_{\xi}=i^*E_{\xi}=E_{\xi}$. 
\par \rm As a component of
$B_{\xi}$,
$E_{\xi}$ has a natural 
structure of ($g-2$)-dimensional scheme: 
\bigskip \it {\bf PROPOSITION 4.4} Assume $\xi$
is general. Then $E_{\xi}$
is reduced and  the intersection $E_{\xi} \cap -E_{\xi}$
is proper. Hence $B_{\xi}$ is reduced.
\par \rm
\it Proof \rm   Let $Z$ be the Zariski closure of 
$B_{\xi}-E_{\xi}$. 
By the previous lemma we have shown 
$-Z \subseteq E_{\xi}$.
Hence
the intersection $Z \cap -Z$ is proper. 
We show that $-Z$ is
reduced and and equal to $E_{\xi}$.
Let us consider the surface
$S_l = \lbrace l+p-x-y, \quad x+y \in C(2) \rbrace$. 
In the next lemma (4.9) it is
shown that the set   
$$
S_l \cap Z
$$  
consists of exactly $2g(g-1)$ distinct points,
provided $\xi$ and $l$ are general.
On the other hand, computing intersection numbers, 
we obtain $(S_l,Z)+(S_l,-Z)$ $=$
$2(S_l,Z)$ $\leq (S_l,B_{\xi}) = 4g(g-1)$,
that is $(S_l,Z) \leq 2g(g-1)$.
Hence $(S_l,Z) = 2g(g-1)$
and  $Z$ is reduced. \par  To show $-Z = E_{\xi}$
observe that
$-Z \subseteq E_{\xi}$, and moreover 
that $(S_l,-Z) = (S_l,E_{\xi})$. Then
$-Z = E_{\xi}$, because $S_l$ is a positive cycle.
\bigskip 
\bigskip \it {\bf COROLLARY 4.5} $E_{\xi}$
has cohomology class $[2\Theta^2]$.
\par \rm We want to explain the
geometrical meaning of the number
$$ 2g(g-1) = (E_{\xi},S_l) = 2[\frac {\Theta^g}{(g-2)!}] 
\tag 4.6
$$ appearing in the proof of the proposition.
With the usual notations we fix a
general
$$ W \in G_l,
$$  where $G_l$ is the 
Grassmannian $Grass(3,H^0(M_l))$ and $M_l = \omega_C(2l+p)$.
Let
$$ f_W: C \to \bold P^2 = \bold PW^* \tag 4.8
$$ be the map defined by $W$. Since $M_l$ is very ample,
we  can assume that: \par
\noindent - $\mid W \mid$ is base-point-free, 
\par
\noindent   - $f_W: C \to f_W(C)$ is birational and
the singular points of $f_W(C)$
are ordinary nodes.\par \noindent 
This follows from generic projection lemma. Of
course we can also assume
$$ W = g_l(\xi)
$$  for some $\xi \in X_l$, (same notations of section 2).
Computing  the number
$\delta$ of nodes of $f_W(C)$, we obtain $\delta = 2g(g-1)$.
\bigskip \it {\bf LEMMA 4.9} 
Let $\xi$ be as above and let $Z$ be the Zariski closure of
$B_{\xi}-E_{\xi}$. There exists a natural bijection
$$ 
b^-: Sing f_W(C) \to S_l \cap Z, \tag 4.10
$$ where $S_l$ is the surface
$S_l = \lbrace l+p-x-y, \quad x+y \in C(2) \rbrace
\subset Pic^0(C)$.
\par \rm
\it Proof \rm  Let $Sing f_W(C) = \lbrace o_1, \dots, o_{\delta} \rbrace$,
($\delta = 2g(g-1)$). For each $o_i$ we consider the two branches 
$$
\lbrace x_i,y_i \rbrace = f_W^{-1}(o_i), 
$$ and the point $e_i = l+p-x_i-y_i \in S_l$. Let us show that
$$ S_l \cap Z = \lbrace e_1, \dots, e_{\delta} \rbrace.
$$ For $l$ general we have $S_l \cap Z \cap E_{\xi} = \emptyset$.
Therefore we can
assume
$h^0(\xi(e)) = 1$, for each $e \in Z \cap S_l$. 
Moreover we know that a point
$e$ satisfying $h^0(\xi(e)) = 1$ belongs to $B_{\xi}$ iff
$h^0(\xi(e-p)) = 1$. Therefore, for a point $e = l+p-x-y \in S_l$,
it follows
$$ e \in Z \Longleftrightarrow h^0(\xi(l-x-y)) = 1.
$$
It is easy to see that the latter condition is 
satisfied if and only if $f_W(x) =
f_W(y)$. That is iff $o = f_W(x) = f_W(y)$ is a node of $f_W(C)$. Then
$$
Z \cap S_l = \lbrace e_1 \dots e_{\delta} \rbrace.
$$ 
This defines a bijective map $b^-:Sing f_W(C) \to B^1_{\xi} \cap S_l$
sending
$o$ to $l+p-f_W^*(o)$. \bigskip 
In the proof of proposition (4.4) we have shown
that $Z=-E_{\xi}$ if $\xi$ is general.
Therefore propositions (4.4) and (4.9) imply
the following
\bigskip \it {\bf PROPOSITION 4.11} 
Assume the pair $(\xi,l) \in X \times C$ is general. Let
$W = g_l(\xi)$ and let $S_l = l+p-C(2)$. 
Then there is a bijection
$$
b^-: Sing f_W(C) \to (-E_{\xi}) \cap S_l
$$  sending $o \in Sing f_W(C)$ to $b^-(o) = l+p-f_W^*o$. 
\par \rm We want to complete the picture by a characterization of
$$ E_{\xi} \cap S_l.
$$ We will only sketch this, 
leaving the proofs as an exercise. \par
For a general pair
$(\xi,l)$ one has $h^0(\xi(l)) = 5$. We consider the ruled surface
$$  R = Proj(\xi)
$$  and the map
$$  g: R \to \bold P^4 = \bold PH^0(\xi(l+p))^*,
$$  which is induced by the evaluation 
$H^0(\xi(l+p))\otimes \Cal O_C \to \xi(l+p)$.
We assume that
$$ g: R \to g(R)
$$ is a birational morphism and that $Sing f_W(R)$
consists of finitely many apparent
double points, (that is points $o$ such that $g^*o$
is reduced of length two). One
can show that a general pair $(\xi,l)$ satisfies this assumption.
Again, the formula
for the number $\delta$ of apparent 
double points says  $\delta = 2g(g-1)$.
Let
$\pi: R \to C$ be the natural projection,
for each $o \in Sing (g(R))$ we have a point
$$ e = l+p-x-y, \quad \text {where $\pi_*g^*o = x+y \in C(2)$.}
$$ It is easy to check that $h^0(\xi(e)) \geq 2$,
therefore $e \in E_{\xi} \cap S_l$.
\bigskip \it {PROPOSITION 4.12} Assume the pair 
$(\xi,l) \in X \times C$ is general. Let \par
\noindent
$g: R \to \bold P^4$ be as above and let $S_l = l+p-C(2)$.
Then   there exists a
natural bijection
$$  b^+: Sing f_W(C) \to E_{\xi} \cap S_l
$$  sending $o \in Sing f_W(C)$ to $b^+(o) = l+p-\pi_*g*o$. 
\par \rm Finally we describe how $E_{\xi}$ parametrizes 
the curves $l-C$ which
are contained in a divisor $D_{\lambda} \in P_{\xi}$.
To a point $e \in E_{\xi}$ we
associate the curve
$$ C_e = \lbrace e+p-x, \quad x \in C \rbrace,
$$ passing through $e$. $C_e$ is contained in at
least one $D_{\lambda}$ of the pencil
$P_{\xi}$. Indeed, since $h^0(\xi(-e-p))\geq 1$, 
there exists at least one exact
sequence
$$
\CD 0 @>>>J{\eta_{\lambda}} @>>>
{\xi} @>{\lambda}>> {\bold C_p} @>>> 0 \\
\endCD
$$ such that $h^0(\eta_{\lambda}(-e-p)) \geq 1$. 
For such a $\lambda$ we have $C_e
\subseteq D_{\lambda}$. Let
$$ F_{\xi} = \lbrace e \in Pic^0(C)
/ e+p-C \subset D_{\lambda} \text {, for some
$\lambda$} \rbrace,
$$
$F_{\xi}$ is a closed set. From the previous remarks it follows 
$$ E_{\xi} \subseteq F_{\xi}.
$$ Assume the intersection scheme 
$$
\Delta_e = C_e \cdot B_{\xi}
$$ is smooth and finite. Then, since
$P_{\xi}$ is a pencil of $2\Theta$-divisors 
and $(\Theta,C_e) = g$,  $C_e \cap
B_{\xi}$ consists of exactly $2g$ points. 
The distribution of these points is the following:
$2g-1$ points on $-E_{\xi}$, the point $e$ on $E_{\xi}$,
(see the next
proposition).\par \noindent Let $e \in E_{\xi}$, 
assume that the determinant map $v:
\wedge^2 H^0(\xi(e)) \to H^0(\omega_C(2e))$
is injective and that $h^0(\xi(e)) = 2$.
As in lemma (3.13), the image of $v$ defines the  divisor
$$ d_e \in \mid \omega_C(2e) \mid.
$$
One can easily check that, as a divisor on $C$,
$$
\Delta_e = d_e+p.
$$
\bigskip \it {\bf PROPOSITION 4.13}  $E_{\xi} = F_{\xi}$
if $C$ is not hyperelliptic.
\par \rm
\it Proof \rm   Assume $e \in F_{\xi}$.
It suffices to show that, for some $\lambda
\in \bold P^1$, one has $h^0(\eta_{\lambda}(-e-p)) \geq 1$.
This implies
$h^0(\xi(-e-p)) \geq 1$ and hence $e \in E_{\xi}$.
Let $Y$ be an irreducible
component of $F_{\xi}$, consider the closed set
$$
\tilde Y = \lbrace (e,\lambda)
\in Y \times \bold P^1/ h^0(\eta_{\lambda}(-e-p))
\geq 1 \rbrace
$$ and the complement $U$ of its projection in $Y$.
We must show that $U$ is empty.
Assume $e \in U$. Then, for some $\lambda = o$, we have:
\par \noindent (i)
$h^0(\eta_o(e+p-x))
\geq 1$,$\forall x \in C$, (ii) $h^0(\eta_o(-e-p)) = 0$. 
\par \noindent Note that
(ii) is equivalent to $h^0(\eta_o(e+p)) = 2$. 
Conditions (i) and (ii) imply that
$\eta_o(e+p)$ contains a subbundle $L \in W^1_d$,
where $W^1_d$ is the
Brill-Noether locus. Clearly, the same holds for $\xi(e+p)$.
Consider the set of pairs
$$ T = \lbrace (e,L) 
\in Y \times W^1_d / h^0(\xi(e+p) \otimes L^*)) \geq 1
\rbrace,
$$ for the dimension of $T$ we have $dim T \geq dim U = dim Y$. 
The image of $T$ under
the map $(e,L) \to L(-e-p)$ is a family $Z$ of line subbundles 
of $\xi$ having degree
$d-1$. Since the fibre of this map is isomorphic to $W^1_d$, it follows 
$$ dim Z \geq dim Y - dim W^1_d.
$$ On the other hand, by proposition (3.10), we have 
$$ dim Z \leq g-(d-1)-1 = g-d.
$$ 
Observe that $F_{\xi} = \cup H_{\lambda}$, where $H_{\lambda} = \lbrace
e \in Pic^0(C) / e+p-C \subset D_{\lambda} \rbrace$. 
It is known that $dim H_{\lambda}
\geq g-3$, moreover the dimension of 
$\lbrace e \in E_{\xi} / e+p-C \subset E_{\xi}
\rbrace$ is $\leq g-3$, (cfr. [BV]
proposition (5.10) and lemma (5.14)). Therefore it
follows 
$$ dim Y \geq g-2.
$$ Assume $C$ is not hyperelliptic.
Then, by Martens theorem, $dim W^1_d
\leq d-3$. This implies $g-d \geq dim Z \geq dim Y -
dim W^1_d \geq g-d+1$, which is
a contradiction. Therefore $U$ is empty and $F_{\xi}
= E_{\xi}$. \par \noindent 
\bigskip 
\bigskip \it {\bf REMARK} \rm $E_{\xi}$ defines a divisor $D_{\xi}$ 
which is naturally
associated to $\xi$. Let $ \Cal C = \lbrace (e,f)
\in E_{\xi} \times Pic^0(C) / f \in
C_e \rbrace$, then 
$$ 
D_{\xi} = \pi_* \Cal C
$$
where $\pi: E_{\xi} \times Pic^0(C) \to E_{\xi}$
is the first projection. 
\par \rm
Let $G = Grass(2,H^0(\Cal O_J(2\Theta))$ 
and let $X_C$ be the moduli space of stable
rank two vector bundles of 
determinant $\omega_C(p)$, 
where $p \in C$. We can
consider the morphism
$f: X_C \to G$ such that $f(\xi) = 
\phi_p(\xi)$ if $det \xi
= \omega_C(p)$. It is natural to 
end this section with the following
\bigskip \it {PROBLEM} \rm Is it true 
that $f(X_C)$ is an irreducible component of the
variety of pencils $P \in G$ having 
reducible base locus?
\par \rm
\centerline{\bf References}
\hfill \par [A-C-G-H] E.Arbarello,M.Cornalba,
P.A.Griffiths,J.Harris,{\it Geometry of
Algebraic Curves,I}, Springer, Berlin (1985).
\hfill \par [B-L]JC.Birkenhake, H.Lange, 
{\it Complex Abelian Varieties}, Springer,
Berlin (1992).
\hfill \par [B1] A.Beauville,
{\it Fibres de rang 2 sur une courbe,
fibr\'e determinant et functions theta}, 
Bull. Soc.Math.France, {\bf
116}(1988),431--448.
\hfill \par [B2] A.Beauville,
{\it Fibres de rang 2 sur une courbe, fibr\'e determinant
et functions theta, II}, Bull. 
Soc.Math.France, {\bf 119}(1991),259--291.
\hfill \par [B3] A.Beauville, 
{\it Vector bundles 
on curves and
generalized theta
functions: recent results and open problems}
in {\it Current topics in complex
algebraic geometry}, Cambridge University 
Press, Cambridge (1995), 17--33
\hfill \par [B-N-R] A.Beauville,M.S.Narasimhan,S.Ramanan, 
{\it Spectral curves and
the generalised theta divisor}, J.reine angew. Math. 
{\bf 398}(1989), 169--179.
\hfill \par [Br] S.Brivio, {\it On rank
2 semistable vector bundles over an
irreducible nodal curve of genus 2},
to appear on Boll. U.M.I.
\hfill \par [B-V] S.Brivio, A.Verra, 
{\it The theta divisor of ${{\Cal SU}_C(2)}^s$
is very ample if $C$ is not hyperelliptic},
Duke Math. J., {\bf 82}(1996), 503--552
\hfill \par [D-N] I.M.Drezet,M.S.Narasimhan, 
{\it Groupes de Picard des vari\'et\'es
des modules des fibr\'es semistable sur les courbes algebriques},
Invent.Math. {\bf
97}(1989), 53--94.
\hfill \par [D-R] U.V.Desale,S.Ramanan,
{\it Classification of vector bundles of rank
two on hyperelliptic curves}, Invent.Math.,
{\bf 38}(1976),161--185.
\hfill \par [L-N] H.Lange, M.S.Narasimhan,
{\it Maximal subbundles of rank two
vector bundles}, Math. Annalen,
{\bf 266}(1983), 55--73
\hfill \par [L1] Y.Laszlo,
{\it A propos de l'espace des modules de fibres de rang 2
sur une courbe}, Math.Annalen, {\bf 299} (1994) 597-608
\hfill \par [L2] Y.Laszlo,
{\it Un th\'eor\'eme de Riemann pur le diviseur Theta
g\'en\'eralis\'e sur 
les espace de modules de fibr\'es stables sur une courbe}, Duke
Math. J., {\bf 64} (1991) 333-347 
\hfill \par [N-R1] M.S.Narasimhan, S.Ramanan,
{\it Moduli of vector bundles on
a compact Riemann surface}, Ann.Math. {\bf 89}(1969), 19--51.
\hfill \par [N-R2] M.S.Narasimhan, S.Ramanan,
{\it $2 \theta $-linear systems on
Abelian Varieties}, Vector Bundles on
Algebraic varieties, Oxford (1989),415--427.
\hfill \par [R] M. Raynaud, {\it Sections des
fibr\'es vectoriels sur une courbe},
Bull.Soc.math.\par \noindent France, {\bf 110}(1982), 103--125.
\hfill \par [S] C.S.Seshadri, 
{\it Fibr\'es vectoriels sur les courbes
alg\'ebriques}, Ast\'erisque, {\bf 96}, (1982) 
\hfill \par
\hfill \par \bf Authors'addresses: \rm
\hfill \par \ S.Brivio, Dipartimento di Matematica,
Universita' di Torino, via Carlo Alberto,10 - 10123 Torino (Italy).
\hfill \par  A.Verra, Dipartimento di Matematica, Universita' di Roma Tre,
largo S. Leonardo Murialdo 1 - 00146 Roma (Italy).
\enddocument